\documentclass[onefignum,onetabnum]{siamonline190516}

\usepackage{lipsum}
\usepackage{amsfonts}
\usepackage{graphicx}
\usepackage{subcaption}
\usepackage{epstopdf}
\usepackage{epigraph}
\usepackage{enumitem}

\usepackage{algorithmic}
\ifpdf
  \DeclareGraphicsExtensions{.eps,.pdf,.png,.jpg}
\else
  \DeclareGraphicsExtensions{.eps}
\fi

\usepackage{enumitem}
\setlist[enumerate]{leftmargin=.5in}
\setlist[itemize]{leftmargin=.5in}

\newsiamremark{remark}{Remark}
\newsiamremark{hypothesis}{Hypothesis}
\crefname{hypothesis}{Hypothesis}{Hypotheses}
\newsiamthm{claim}{Claim}

\newcommand{\mA}{\mathbf{A}}

\newcommand{\transpose}     {^{\mbox{\scriptsize \sf T}}}
\newcommand{\mB}{\mathbf{B}}
\newcommand{\mC}{\mathbf{C}}
\newcommand{\mY}{\mathbf{Y}}
\newcommand{\mM}{\mathbf{M}}
\newcommand{\mX}{\mathbf{X}}

\newcommand{\mR}{\mathbf{R}}
\newcommand{\mP}{\mathbf{P}}
\newcommand{\mF}{\mathbf{F}}

\newcommand{\mT}{\mathbf{T}}
\newcommand{\mL}{\mathbf{L}}
\newcommand{\mS}{\mathbf{S}}
\newcommand{\dnnz}{\mathit{nnz}}

\usepackage{mathtools}

\usepackage{listings}

\renewcommand{\matrix}[1]{{\bf #1}}

\headers{The Ubiquitous Sparse Matrix-Matrix Products}{Aydin Bulu\c{c}}

\title{The Ubiquitous Sparse Matrix-Matrix Products}

\author{Ayd\i n Bulu\c{c}\thanks{Lawrence Berkeley National Laboratory and University of California, Berkeley 
  (\email{abuluc@lbl.gov}, \url{https://people.eecs.berkeley.edu/\~aydin}).}}

\usepackage{amsopn}

\ifpdf
\hypersetup{
  pdftitle={The Ubiquitous Sparse Matrix-Matrix Products},
  pdfauthor={Ayd\i n Bulu\c{c}}
}
\fi

\begin{document}

\maketitle

\begin{abstract}
  Multiplication of a sparse matrix with another (dense or sparse) matrix is a fundamental operation that
  captures the computational patterns of many data science applications, including but 
  not limited to graph algorithms, sparsely connected neural networks, graph neural networks, 
  clustering, and many-to-many comparisons of biological sequencing data. In many
  application scenarios, the matrix multiplication takes places on an arbitrary algebraic semiring where 
  the scalar operations are overloaded with user-defined functions with certain properties or a more
  general heterogenous algebra where even the domains of the input matrices can be different. Here, we
  provide a unifying treatment of the sparse matrix-matrix operation and its rich application space including 
  machine learning, computational biology and chemistry, graph algorithms, and scientific computing.
  \end{abstract}


\epigraph{``{\it I am called Ubik, but that is not my name. I am. I shall always be."} }{\textit{Philip K.\ Dick}}

\section{Introduction}
Consider the operation of multiplying a sparse matrix with another matrix. The second matrix can
be sparse or dense, and this often changes both the applications of this ``sparse matmul'' operation and the best
algorithms to compute the product. Mathematically, however, the more important distinction is whether the operation is truly mapping a list of vectors from one
a vector space to another. This is often the case when the second matrix is ``tall-and-skinny'', meaning that it has many more rows than columns.

We call a matrix sparse if it has sufficiently many ``null'' elements such that storing
and treating it as a dense matrix would be computationally wasteful or infeasible. The null element for classical numerical linear algebra is
often the real number $0$, however it will be different in other application domains. Not only the null element can be different, but
also the whole algebra the product is computed on can be different. This algebra can in principle be anything where the computation is well defined, 
but it often is referred to as a semiring $\mathbb{S}$, which is a ring without the additive inverses. The added constraint of computing the ``sparse matmul'' when
additive inverses are not available rule out many ideas from fast matrix multiplication such as Strassen's algorithm~\cite{strassen1969gaussian} and those based on arithmetic progression~\cite{coppersmith1987matrix}.

The sparse matmul operation is then $\mY \gets \mA \mX$ where $\mA \in \mathbb{S}^{M \times K}$ is a sparse matrix but $\mX \in \mathbb{S}^{K \times N}$ and $\mY \in \mathbb{S}^{M \times N}$ can be either sparse or dense. To allow extensions to tensor multiplication, we refer to the inner dimension $k$ as the \emph{contraction dimension}. 
The \emph{aspect ratio} of a $M\times N$ matrix is $M/N$, which is the ratio of its number of rows to its number of columns. 
When $\mX$ is also sparse, the sparse matmul operation is commonly referred to as SpGEMM in the literature, regardless of the aspect ratio of the matrix. 
In this case when both $\mA$ and $\mX$ are sparse, the output $\mY$ can be represented as a dense or sparse matrix, depending on its nonzero density.
When $\mX$ and $\mY$ are dense but have a high aspect ratio, the operation becomes a multi-vector generalization of the sparse matrix times dense vector product, as is
commonly referred to as SpMM in the literature. The naming of the case where $\mX$ and $\mY$ are dense but square is not well established but a few publications
refer to it as the SpDMMM or SpDM$^3$ in short. 

SpGEMM and SpMM are part of the upcoming Sparse BLAS standard~\cite{abdelfattah2024interface}. They can also be invoked internally within the \texttt{GrB\_mxm()} function of the GraphBLAS standard depending on the sparsity of the input~\cite{brock2021graphblas}. The GraphBLAS interface does not specify the sparsity of its matrices and vectors, hence the implementations are free to choose a sparse or dense representation for a matrix; and use different internal functions accordingly. 

A function that is closely related to SpMM is the sampled dense-dense matrix multiplication, also known as SDDMM. This function multiplies two dense matrices but only for the set of nonzero indices required by a dense sampling matrix. Concretely, this operation can be written as $\mY \langle\mM\rangle \gets \mA \mX$ where $\mA$ and  $\mX$ are dense matrices, and the sparse sampling matrix $\mM$ denotes which indices of the output needs to be computed. It is mathematically equivalent to $\mY \gets \mM \odot \mA \mX$ where $\odot$ is the Hadamard product. The GraphBLAS standard calls the sampling matrix a ``mask'' and generalizes its use such that it can be applied to cases beyond SDDMM, including a masked version of SpGEMM.

For clarity, let us enumerate all the special cases of what is often vaguely called ``sparse matrix multiplication'' in the literature following the $\mY \gets \mM \odot \mA \mX$ notation:
\begin{enumerate}[leftmargin=10pt]
\item \textbf{SpGEMM}: $\mA$, $\mX$ are sparse, $\mY$ can be sparse or dense (depending on shape/density)
\item \textbf{Masked-SpGEMM}: Same as SpGEMM, with a mask ($\mM$) present.
\item \textbf{SpMM}: $\mA$ sparse, $\mX$ and $\mY$ dense and tall-skinny, often no mask.
\item \textbf{SDDMM}: $\mA$, $\mX$ are dense, $\mM$ present, $\mY$ sparse due to mask.
\item \textbf{SpDM$^3$}: $\mA$ sparse, $\mX$ and $\mY$ dense and square, often no mask.
\item \textbf{SpMV}: degenerate case of SpMM with $\mX$ and $\mY$ having 1 column.
\item \textbf{SpMSpV}: degenerate case of SpGEMM with $\mX$, $\mY$ (and possibly $\mM$) having 1 column.
\end{enumerate}

 \begin{figure}[thb]
  \centering
 \includegraphics[width=0.9\linewidth]{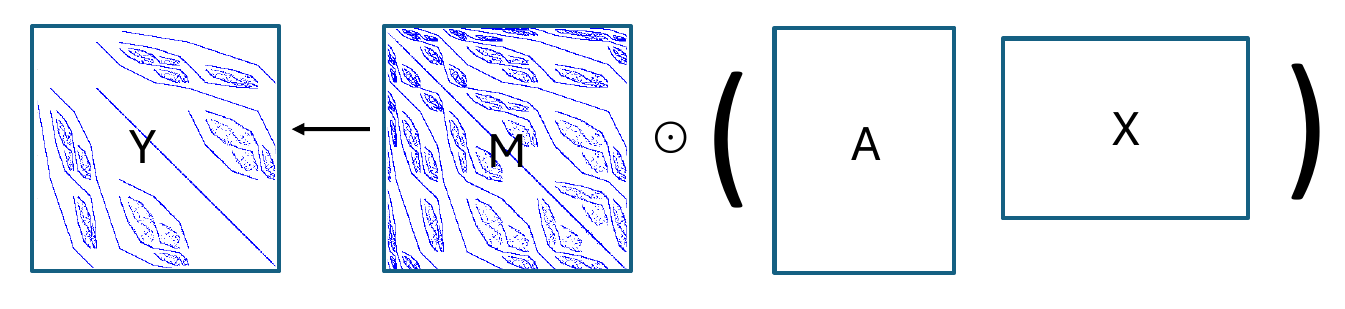}
 \caption{An example for (masked) sparse matrix multiplication. In general, $\mA$ and $\mX$ can each be dense or sparse, but in this example at least one is sparse because the sparsity pattern of $\mY$ is a strict subset of the sparsity structure of $\mM$. }
 \label{fig:sparsemult}
 \end{figure}
 
Notice that if the mask is present, it is always sparse. Figure~\ref{fig:sparsemult} gives an illustrative example. We will not cover SpMV and SpMSV in this paper but they are of utmost importance to scientific computing (especially iterative solvers) and graph algorithms (especially graph traversals). We plan to cover SpMV and SpMSV in a follow-up paper.

This paper presents the use cases and application scenarios of matrix-matrix multiplication where at least one input matrix is sparse. 
This area has been subject to intensive and accelerating research in the last decade. We will not be surveying algorithmic and architectural research here. 
Many papers that develop new algorithms and architectures for sparse matmul often copy motivations from prior work without verifying the exact use cases. As a result, they often fail 
to utilize realistic benchmarks in their experimental evaluation. By identifying the sources of matrices and carefully 
defining the basic features such as the sparsity and matrix sizes, our paper will help algorithm designers as well as computer architectures who develop specialized hardware for such problems~\cite{pal2018outerspace,zhang2020sparch}. The majority of the use cases described in this article are published in the literature but we also present a few new use cases.

In addition to providing realistic use cases for sparse matrix developers, this work has an equally important complementary goal. We provide application developers a new way to look at certain patterns
of computation in the lens of sparse matrices. This lens provided three benefits in the past: (1) it opened up new optimizations by thinking in the matrix level, (2) it allowed novel parallelization techniques for problems where the parallelism opportunity was not obvious, and (3) it enabled applications to use optimized sparse matrix codes. 

\section{Performance Considerations}

Under certain idealized assumptions, the operational intensities (OI) of SpMM, SpGEMM, and SpDM$^3$ can be compared for estimating best case performance differences. For instance, assume that we have a square $N$-by-$N$ sparse matrix $\mA$ has $\dnnz(A)=Nd$ nonzeros that are uniformly distributed and we would like to multiply this matrix with a dense tall-skinny matrix. If the short dimension of the dense matrices are $\Omega(d)$, then 
$$ \text{OI(SpMM)})= \frac{2 N d \, \Omega(d)}{dN + 2N \, \Omega(d)} \approx O(d).$$

This above bound for SpMM also applies to SpDM$^3$  as can be seen easily by plugging $n$ for the second dimension of the dense matrices. Also, because SDDMM and SpMM are shown to be algorithmically identical~\cite{bharadwaj2022distributed} (i.e., any algorithm for one can be mechanically converted to an algorithm for the other with the same costs), the above bound for SpMM also applies to SDDMM. For SpGEMM, the best case OI is often measured using \emph{compression ratio}, which is the ratio of required ``sparse flops'' to the number of nonzeros in the output. Sparse flops is the number of nontrivial arithmetic operations required for computing SpGEMM, i.e., $a_{ij} b_{jk}$ products for which both $a_{ij}$ and $b_{jk}$ are nonzero. For the compression ratio to meaningfully approximate OI, the output sparse matrix needs to have at least as many nonzeros as the input sparse matrices, which often is the case~\footnote{AMG restriction and prolongation is a notable exception, which is covered in Section~\ref{sec:spgemmsci}}. 

For SpGEMM, the compression ratios often range from 1 to 32~\cite{nagasaka2019performance}. For SpMM, the short dimension of the dense matrices ($K$) often range from 8 to 512 ($16-512$ for graph neural networks~\cite{huang2020ge}, $8-96$ for block iterative solvers~\cite{aktulga2016high}). Hence, there is an order of magnitude difference in OI between SpGEMM and SpMM in favor of SpMM. We emphasize that these back-of-the-envelope OI calculations are only relevant for single-node execution. When running on distributed-memory, outer-level parallelization techniques such as 1.5D, 2D, and 3D algorithms can significantly alter the OI of node-level SpMM and SpGEMM kernels~\cite{selvitopi2021distributed}.  There are additional issues that make SpGEMM harder to extract optimal performance compared to SpMM/SpDM$^3$: the control flow of SpGEMM is often more complex, its maintains more complex internal data structures, and the output size is not known in advance. 

We verify our conjecture that SpMM has higher OI and easier to get good performance via experimental results presented in two recent papers. Both papers used 500+ matrices from the SuiteSparse collection. SpMM on NVIDIA P100 achieved a maximum of over $1000$ GFlops when $K=512$, with the median matrix achieving $\approx 400$ GFlops~\cite{jiang2020novel}. The same P100 could only achieve a maximum of $100$ GFlops with the median matrix achieving only $3-4$ GFlops, even with careful ML guided parameter tuning over 6 different implementations~\cite{wei2024predicting}. 

\section{The SpGEMM case}
\label{sec:spgemm}

Computing the product of two sparse matrices, which itself results in a sparse output matrix, has enjoyed early applications in
theoretical computer science~\cite{kaplan2006colored}. Fundamentally, there are three distinct cases of SpGEMM worth mentioning. 
\begin{enumerate}
\item The canonical case involves computing the unrestricted $\mY \gets \mA \mX$ when all three matrices are stored in sparse format.
\item When the output $\mY$ is expected to be above certain density, it is often advantageous to store it in a dense format, as this simplifies the computation.
\item A significant number of graph algorithms and certain machine learning algorithms require only a subset of the output indices. We call this the Masked-SpGEMM and use
the notation $\mY \langle\mM\rangle \gets \mA \mX$ where ``mask'' $\mM$ is a sparse matrix denoting which indices of the output needs to be computed. 
\end{enumerate}

\subsection{SpGEMM for Scientific Computing}
\label{sec:spgemmsci}
\begin{figure}[thb]
  \centering
 \includegraphics[width=0.9\linewidth]{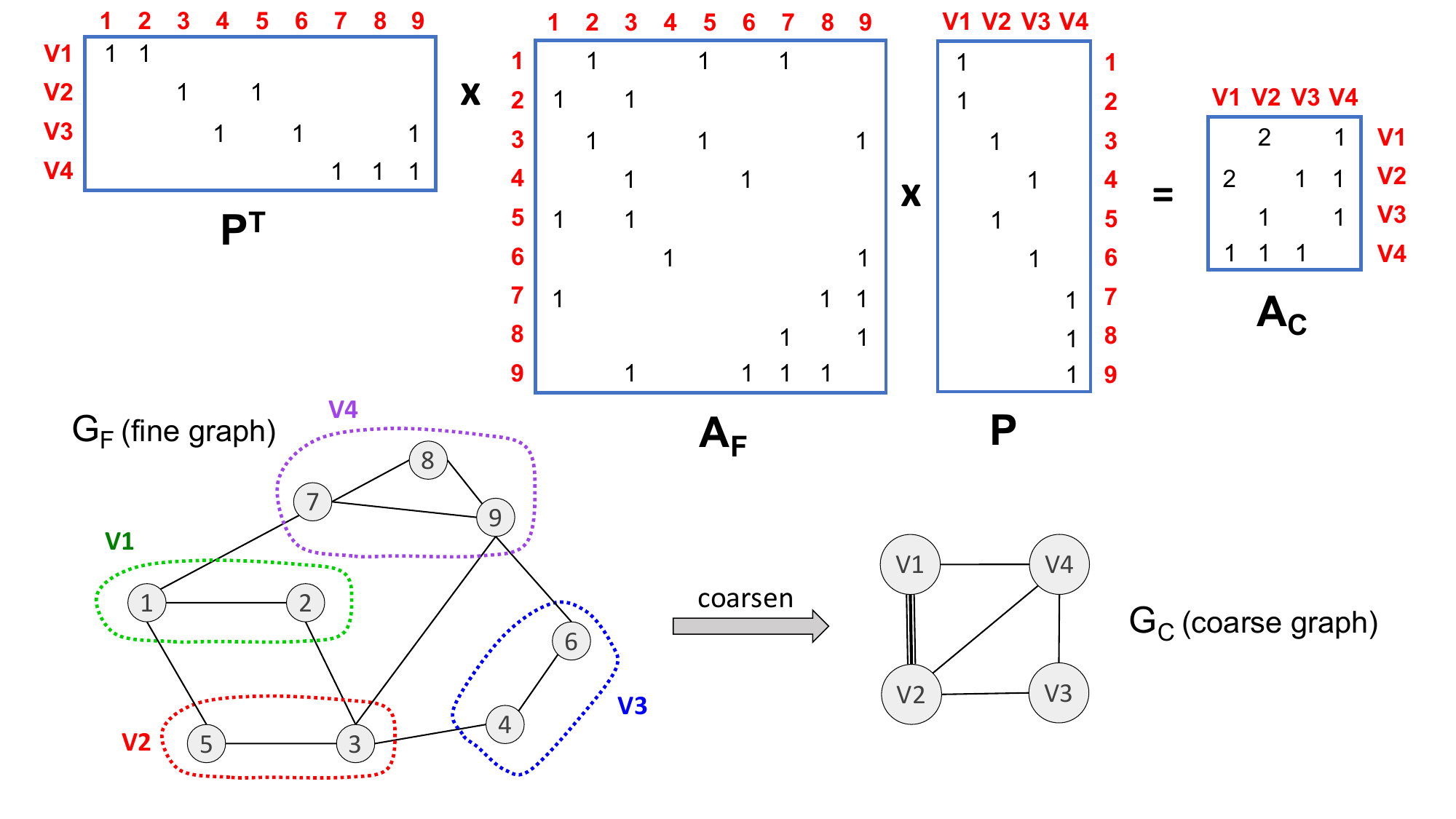}
 \caption{Graph contraction using triple sparse matrix product $\mA_C \gets  \mP\transpose \, \mA_F \, \mP$. The input graph is the fine graph $G_F$, whose adjacency matrix is $\mA_F$. The output graph is the coarse graph $G_C$, whose adjacency matrix is $\mA_C$. Contraction mechanisms other than using the standard $(+, \times)$ field is possible as well as other prolongation matrices.  A fine graph with unit weights is shown for simplicity. Missing locations in the matrices are considered zero.}
 \label{fig:contraction}
 \end{figure}

The earliest practical applications of SpGEMM come from numerical linear algebra, however. There, the triple product $\mC \gets \mR \mA \mP$, also known as the Galerkin product, is
used to compute grid restriction and prolongation (interpolation) during the setup phase of the algebraic multigrid (AMG) method for solving systems of linear equations~\cite{bienz2016reducing}. When 
the system $\mA$ is symmetric, $\mR = \mP\transpose$. 
The bottleneck in AMG, until a decade ago, was the solve phase instead of the setup phase. This has changed recently due to advances in the solve phase,
resulting in a renewed interest in SpGEMM in general and the triple product in particular by researchers working on AMG solvers~\cite{bell2012exposing, ballard2016reducing}.
In this use case of SpGEMM, an unstructured fine grid with $n$ points is converted to another coarser grid with $m$ points using the prolongation matrix $\mP \in \mathbb{R}^{n \times m}$.
The number of fine grid points $m$ depends on the restriction factor of the AMG solver but assuming each fine-grid point only contributes to one and only coarse-grid point, each row of $\mP$ has a single nonzero, resulting in $\dnnz(\mP) = n$. This operation is illustrated in Figure~\ref{fig:contraction}. While this application of SpGEMM was known implicitly by the community for a long time, Bell et al.~\cite{bell2012exposing} gave one of the first
concrete descriptions of an AMG solver that explicitly uses SpGEMM. Their prolongation/restriction operator is based on Maximal Independent Set, which exposes high degree of parallelism. Recently, Li et al.~\cite{li2021new} proposed other prolongation schemes that also use SpGEMM as their work engine. One special characteristic of AMG prolongation/restriction case (and in general graph coarsening)  is that the number of nonzeros in the output $\mC$ (as well as the intermediate products $\mA \mP$ or $\mR \mA$, if the output is computed via two successive SpGEMM calls) can be smaller than the input matrix $\mA$, which effects efficient algorithm design for SpGEMM for these problems~\cite{ballard2016reducing,wang2024optimization,naumov2015amgx}.

Graph contraction is not limited to solving systems of linear equations. For example, it is an indispensable step in multilevel graph partitioning~\cite{hendrickson1995multi}. It is also used in the famous randomized graph cuts
algorithm of Karger and Stein~\cite{karger1996new}. Recently, graph coarsening is used to accelerate various graph learning algorithms~\cite{ying2018hierarchical}.   

\subsection{SpGEMM for Data Analytics}
Computing pairwise distances among a set of data points is one of my most fundamental operations in data analytics. There are many popular distances, such as the Cosine distance, the Manhattan distance, the Chebyshev distance, to name a few. More of than not, the features associated with each data point is sparse and the distance computation can be written mathematically as SpGEMM on a semiring, namely the computation of the \emph{Gram matrix}: $\mA \mA\transpose$ where $\mA$ is the sparse $\lvert \textit{samples} \rvert$-by-$\lvert \textit{features} \rvert$ input matrix.

However, not all distance computations satisfy the same axioms beyond those minimally required for a semiring. Recall that a semiring is composed of two monoids: the multiplicative monoid $\otimes$ and the additive monoid $\oplus$.  Further, $\otimes$ has to distribute over $\oplus$. While additive identity $\mathit{id}_\oplus$ is required for a semiring, the multiplicative one is not. It is often the case that $\mathit{id}_\oplus$ is also the multiplicative annihilator, i.e. $a \otimes \textit{id}_\oplus = 0$ for any $a$ in the semiring $\mathbb{S}$. As Nolet et al.~\cite{nolet2021semiring} demonstrated, the many-to-many distance calculation might require more than the basic SpGEMM if the multiplicative annihilator is different than additive identity. Consider the case where $\textit{id}_\otimes = 0$ as well, which is the semiring for the Manhattan distance. Then, solely iterating over the intersection of nonzero entries from each column is not sufficient to compute pairwise distances correctly because $a \otimes \textit{id}_\otimes = a$. Instead, the implementation needs to iterate over the union of nonzeros. This increases the complexity bounds and renders off-the-shelf row-by-row and column-by-column SpGEMM implementations~\cite{nagasaka2019performance} unsuitable for the task. 

A similar computation is needed when constructing a graph from a hypergraph. Consider the aforementioned samples-by-features matrix $\mA$ as a hypergraph $\mathbf{H(A)}$, then computing the Gram matrix is equivalent to computing the so-called clique expansion~\cite{agarwal2006higher} of this hypergraph. Liu et al.~\cite{liu2022high} considers the generalized case where there is an input parameter $s$, which acts as a threshold of shared hyperedges needed between pairs of nodes to have an edge in the final graph expansion. They name this problem the $s$-line graph expansion of a hypergraph. The larger the $s$ value is, the sparser the final graph gets. They present algorithms that beat the SpGEMM formulation thanks to their ability to only compute one half of the output, due to the symmetry of the output, and their ability shortcut the computation when the nodes share less than $s$ nodes. We note that these are simple modifications that can be used to optimize an open source SpGEMM code for particular applications, which are generally performed in practice (e.g. the symmetry optimization in BELLA~\cite{guidi2021bella} and the minimum number of shared $k$-mers optimization in PASTIS~\cite{selvitopi2020distributed}, applications we will describe in the next section).

Data clustering is another area where SpGEMM has been employed, starting from the year 2000. Arguably, the most famous clustering method that uses SpGEMM is Markov clustering (MCL)~\cite{dongen2000graph}, a graph clustering algorithm
that is based on random walks. MCL has been especially successful in clustering protein-protein interaction~\cite{krogan2006global} and finding protein families using protein-protein similarity networks~\cite{enright2002efficient}. This success can intuitively be attributed to MCL's
ability to utilize information from other vertices that are many hops away. In one MCL iteration, each vertex shares its normalized probability mass with its immediate neighbors using SpGEMM. This ``expansion'' step is immediately followed by inflation and pruning steps that control convergence and sparsity, respectively. Algebraically, the expansion step is simply squaring the sparse matrix that represents the graph. The nonzero structure of the matrix is equivalent to its adjacency matrix but its values are different because the algorithm normalizes columns to make the input column stochastic~\cite{van2008graph}.  
After $k$ iterations of MCL, each vertex can potentially learn from all of its $k$ hop neighbors.  While protein interaction networks are often modest in size, the protein similarity networks can be huge~\cite{schulz2020giant}, necessitating high-performance parallel implementations of the MCL algorithm. One such work is the HipMCL~\cite{azad2018hipmcl}, which is a distributed-memory version of MCL that can cluster trillion edge graphs on large-scale CPU and GPU-equipped clusters~\cite{selvitopi2020optimizing}. In fact, HipMCL was used in a recent large-scale metagenomics study that more than doubled the number of known protein families~\cite{nature23}. 

SimRank~\cite{jeh2002simrank} is another popular application of SpGEMM in data analysis. SimRank measures structural similarity between nodes in a directed graph, based on the principle that two objects are similar if they are referenced by similar objects. The calculation is iterative: in each iteration, the similarity score for each pair of nodes is updated based on the similarities among the in-neighbors of those nodes. This process naturally leads to a formulation involving matrix operations, specifically, repeatedly multiplying sparse matrices representing the graph structure (i.e., the adjacency matrix $\mA$) and similarity relations (the similarity matrix $\mathbf{S}$). Both $\mA$ and $\mS$ are sparse. Concretely, the iterative SimRank formula can be written as 
$$ \mathbf{S} =\max \{ c( \mA^T  \mathbf{S}  \mA ), \mathbf{I}  \},$$

\noindent
where $\mathbf{I}$ is the identity matrix and $c$ is a damping factor.

\subsection{SpGEMM for Computational Biology}
In addition to the MCL algorithm, which is most popular in computational biology, there are several other applications of SpGEMM in analyzing biological data. 
The majority of these applications involve the computation of the \emph{Gram matrix}: $\mA \mA\transpose$, similar to the distance computation problem described in the previous section. 
The Gram matrix itself is sparse when the task is to detect
overlapping pairs of DNA sequences from whole genome sequencing (e.g., BELLA~\cite{guidi2021bella}) or to identify candidate pairs of proteins that are potentially homologous from a large database (e.g., PASTIS~\cite{selvitopi2020distributed}). 
However, Gram matrix can also be dense when
the task is to compute similarity between whole genomes or metagenomes (e.g., GenomeAtScale~\cite{besta2020communication}) because two genomes or metagenomes almost always share some
degree of sequence similarity.  GenomeAtScale can also be used to compute the Jaccard similarity among other sets of data points, beyond just metagenome sequences. 

For all the use cases described in this subsection, the features are $k$-mers, which are subsequences of length $k$. There is often a sampling procedure involved to choose informative $k$-mers. 
The use of SpGEMM can be expanded beyond finding exact sequence matches by introducing a substitute $k$-mer matrix $\mathbf{S}$, allowing approximate matching between protein sequences by computing
$\mA \mathbf{S} \mA\transpose$~\cite{selvitopi2020distributed}. 

Other than these applications that compute sequence similarity and the MCL algorithm, there are other applications of SpGEMM in bioinformatics.
Guidi et al.~\cite{guidi2021parallel} used repeated SpGEMM to perform transitive reduction of overlap graphs to string graphs, which are kinds of graphs used in genome assembly. 
Jain et al.~\cite{jain2019validating} developed a new algorithm to validate the distance constraints for mapping paired-end reads, which are products of a popular sequencing technology, to genome variation graphs. Their
algorithm uses SpGEMM to take advantage of the sparsity in the reference genome variation graphs.

\subsection{SpGEMM for Computational Chemistry}

In modern electronic-structure codes (e.g., CP2K~\cite{kuhne2020cp2k}, ONETEP~\cite{prentice2020onetep}, CONQUEST~\cite{nakata2020large}) the dominant cost in linear-scaling density-functional theory (DFT) and 
related linear-scaling wave-function methods is multiplying two sparse matrices -- the Hamiltonian $\mathbf{H}$ with the current estimate of the density matrix $\mathbf{P}$. 

The Hamiltonian matrix $\mathbf{H}$ represents the quantum mechanical operator of the system expressed in a localized atomic-orbital basis.
The density matrix $\mathbf{P} = f(\mathbf{H})$, encodes the entire electronic structure -- the occupation of quantum states across the system. The density matrix
determines all ground-state properties in DFT, and is the object we want to compute. 
In traditional DFT for small systems, the density matrix $\mathbf{P}$ can be formed by diagonalizing $\mathbf{H}$, finding its eigenvalues/eigenvectors, and 
then occupying the lowest states up to the number of electrons~\cite{finkelstein2023fast}.
For large systems ($N > 10,000$ atoms), dense diagonalization is computationally prohibitive as it scales with $O(N^3)$. 

Linear-scaling, or $O(N)$, methods avoid explicit diagonalization by (1) expanding the function $f(\mathbf{H})$ using polynomials (Chebyshev, Taylor, or purification polynomials), and (2)
iteratively building $\mathbf{P}$ by repeated sparse-sparse multiplications with $\mathbf{H}$, starting from an initial guess $\mathbf{P_0}$. 
In density-matrix-purification and polynomial-expansion algorithms, each application is of the form 
$\mathbf{P_{n+1}} = f(\mathbf{P_n}, \mathbf{H})$ where $f$ uses matrix products to push the eigenvalues of 
$\mathbf{P}$ towards the desired occupancy (0s and 1s). This process drives the matrix $\mathbf{P}$ toward idempotency, as required for the ground-state density matrix. 
Hence, the ground-state density matrix is built without explicit diagonalization~\cite{niklasson2002expansion, schutt2016gpu}.
Note that these transformations manipulate the eigenvalues indirectly -- no eigenpairs are ever explicitly calculated~\cite{kim2016perspective}. 

An important consequence of this algorithm is the repeated multiplication of the same sparsity pattern: 
the same $\mathbf{H}$ is multiplied hundreds of times with a different iterate $\mathbf{P_n}$, which are approximations to the true density matrices. 
This provides opportunities and justification for optimizing the storage of the sparse matrix $\mathbf{H}$, potentially through blocking, reordering, or both.
Similar to the MCL algorithm described earlier, linear-scaling quantum chemistry codes actively truncate elements to avoid blowing up memory.

Quantum chemistry SpGEMM algorithms are specifically designed around block-sparse matrix formats that reflect the atomic orbital basis structure. 
Kohn's nearsightedness principle~\cite{kohn1996density} of quantum mechanical systems 
imply that the interactions decay with distance. 
Therefore, the truncation mechanism can be based on physical distance cutoffs~\cite{challacombe2000general}, 
in addition to the more classical approach of setting any entry smaller than a threshold to zero~\cite{mohr2014daubechies}. 
A class of algorithms called sparse approximate matrix multiply (SpAMM)~\cite{bock2013optimized} are designed to avoid
the temporary blowup of memory that happens after each multiplication but before truncation. SpAMM codes have the the truncation mechanism built into the algorithm, hence never materializing
entries outside the physical distance cutoff.  

The nonzeros in sparse matrices representing particle interactions are clustered together if they are reordered according to a space filling curve~\cite{bock2013optimized}.
Therefore, quantum chemistry SpGEMM algorithms can exploit the fact that quantum chemistry matrices naturally organize into small dense blocks (typically $6\times6$ to $50\times50$) 
corresponding to atom-centered basis functions, rather than treating each matrix element independently. Many codes follow this block-sparse paradigm~\cite{borvstnik2014sparse}. More info on linear-scaling methods can be found in the survey of Bowler and Miyazaki~\cite{bowler2012methods}.

We note that similar to quantum chemistry, block sparse matrices also come into play in deep learning~\cite{child2019generating}. However, in that case the second matrix is often dense, 
resulting in the SpMM primitive we cover in the subsequent sections.

\subsection{SpGEMM for Graph Analytics}
Perhaps the most widely mentioned use case of SpGEMM is in the graph analytics domain. This is partly due to the SpGEMM being one of the most important specializations of the GraphBLAS' {\tt GrB\_mxm} primitive, and partly due to its successful use in scaling several important graph algorithms. The previously mentioned MCL algorithm is potentially a graph analytics code but we chose to present it in the 
biology subsection due to its popularity in that domain. 

\begin{figure}[thb]
 \centering
 \subcaptionbox{1st step}{%
 \includegraphics[width=.48\linewidth]{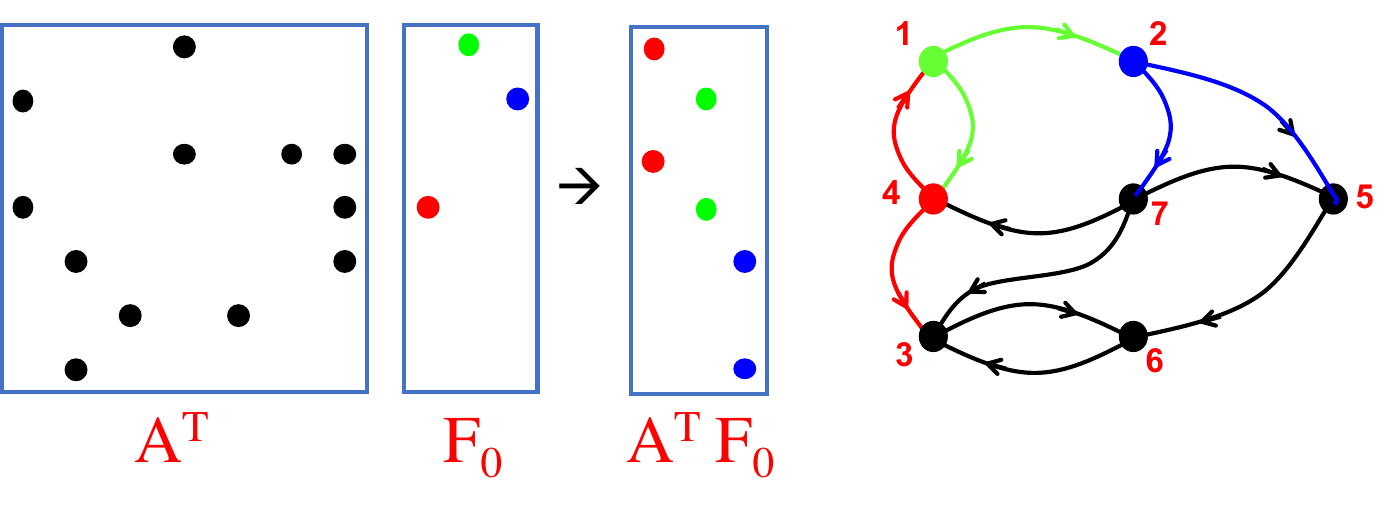}%
 }\hfill
 \subcaptionbox{2nd step}{%
 \includegraphics[width=.48\linewidth]{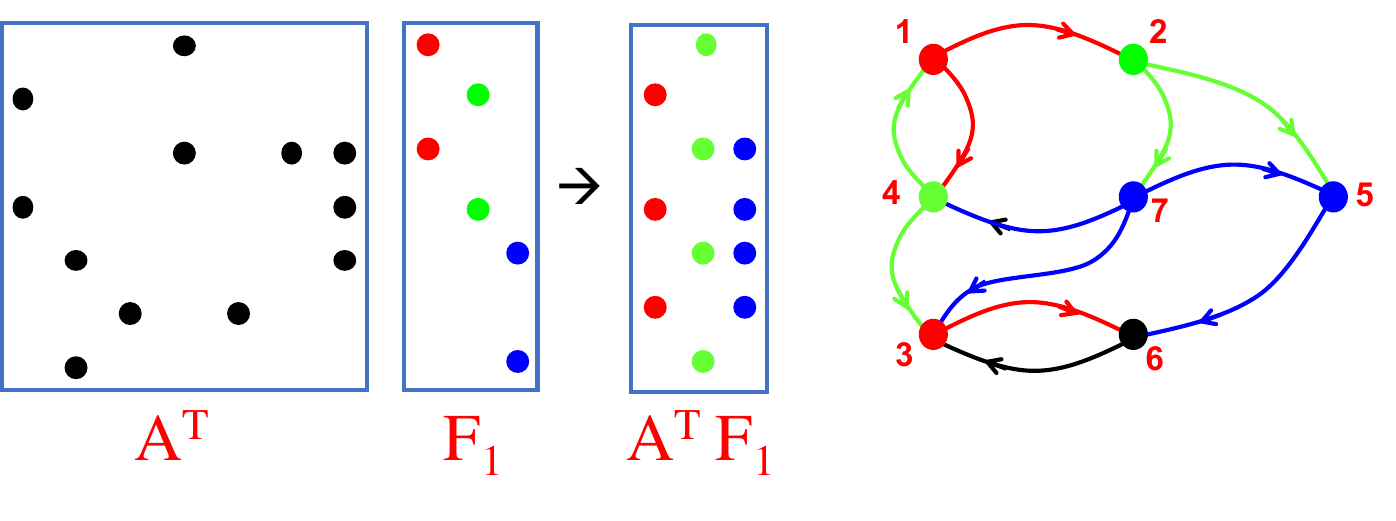}%
 }
 \caption{The first two steps of the multi-source BFS traversal using SpGEMM. $\mA$ is the sparse adjacency matrix of the graph and $\mF_i$ is the set of $i$th frontiers.}
 \label{fig:bc_spgemm}
 \end{figure}

One well-known graph analytical operation that benefits from SpGEMM is the computation of betweenness centrality (BC) on unweighted graphs. BC is a measure of centrality that is based on shortest paths, and is often used to rank entities in a graph. Brandes' algorithm~\cite{brandes2001faster} computes BC efficiently using single-source shortest path (SSSP) computations. The exact BC requires SSSP from every source vertex but using a small fraction of starting vertices often provides a reasonable approximation~\cite{bader2007approximating}. Since breadth-first search (BFS) is sufficient to compute the shortest paths on unweighted graphs, highly-parallel BC implementation perform multi-source BFS (MSBFS). MSBFS maps to SpGEMM with minimal modifications, as illustrated in Figure~\ref{fig:bc_spgemm}. 

Two crucial properties of linear-algebraic BC implementations~\cite{bader2006designing, bulucc2017design} are (1) their use of SpGEMM in a doubly-nested loop, the outer loop on batches of source vertices and the inner loop executing as many times as the diameter of the subgraph induced by the BFS trees sourced at the current batch, and (2) their use of SpGEMM in both the forward MSBFS phase as well as the backwards tallying phase, albeit with different semirings. These two properties in turn make BC performance completely dependent on the SpGEMM performance, as illustrated by large-scale experimental results of Combinatorial BLAS (CombBLAS)~\cite{bulucc2011combinatorial} and Cyclops Tensor Framework (CTF)~\cite{solomonik2017scaling} implementing BC using SpGEMM.  

Other than productivity benefits of using an existing primitive for implementing a complex graph algorithm, this linear algebraic formulation has the unique property of exploiting all available parallelism in BC: (1) multiple independent graph traversals, (2) all vertices in the frontiers of each independent traversal, and (3) all outgoing edges of each one of the frontier vertices. Subsequent papers proposed special methods~\cite{prountzos2013betweenness} for exploiting all three levels of parallelism in BC, highlighting the hidden benefit of the linear-algebraic approach. 

\begin{figure}[thb]
  \centering
 \includegraphics[width=0.9\linewidth]{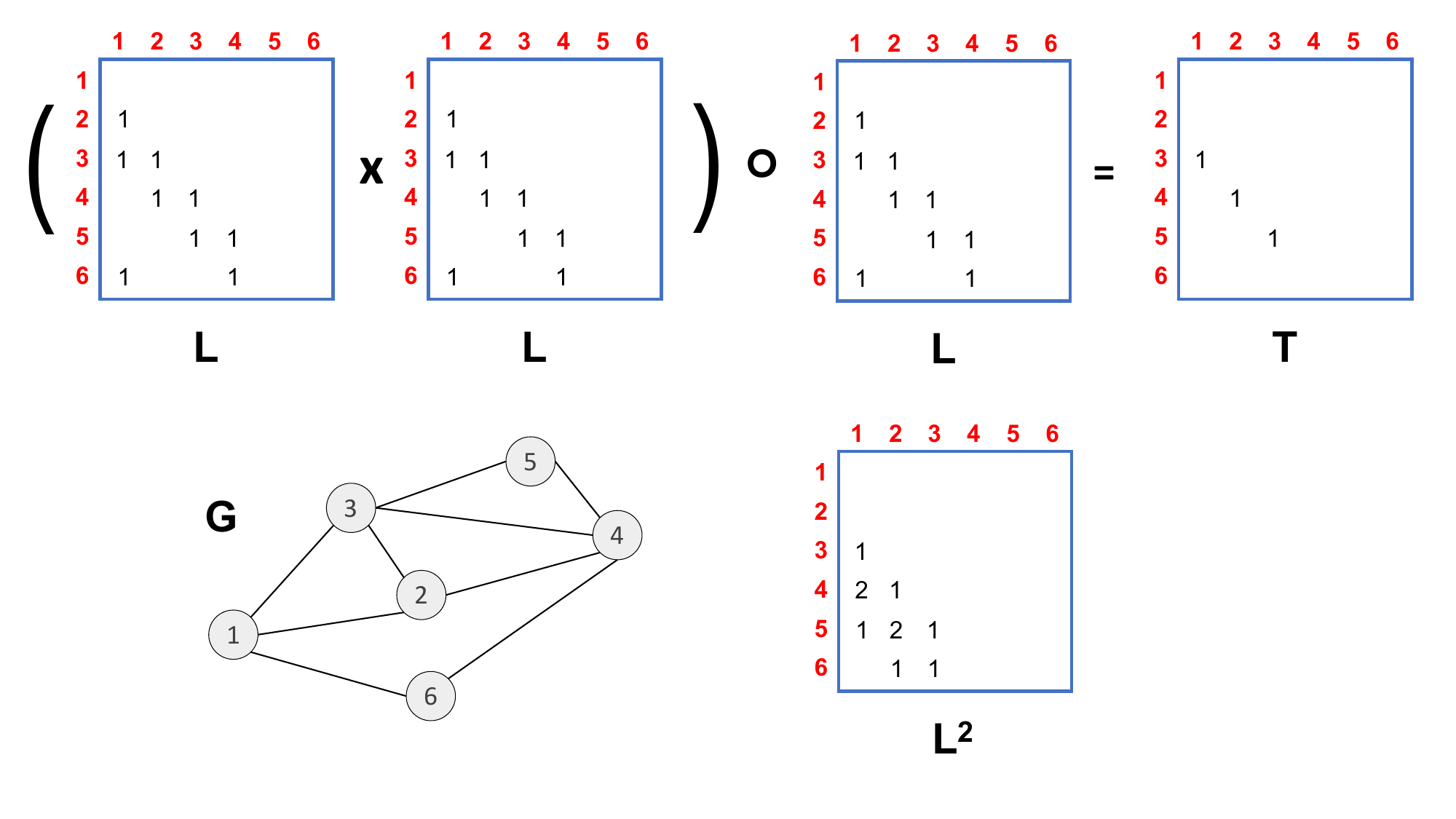}
 \caption{Triangle counting via SpGEMM. Each nonzero $\mL^2(i,j) \neq 0$ shows the number of triangles that could form if $\mL(i,j) \neq 0$. In other words, it lists the wedges around the potential edge $e(i,j)$.
 Computing the Hadamard product of $\mL^2$ with $\mL$ tests whether the triangles actually form (top of the figure).} 
 \label{fig:tc}
 \end{figure}

The linear algebraic formulation of triangle counting, while seemingly a unique algorithm at first sight, is actually another formulation of graph traversal as matrix multiplication. 
Consider perhaps the most efficient formulation~\cite{wolf2017fast} where entries $\mT(i,j) \neq 0$ of $\mT = \mL^2 \circ  \mL$ give the number of triangles centered around that edge $e(i,j)$ where $\circ$ is the Hadamard product. 
The total number of triangles in the whole graph is simply $\sum_{i,j}{\mT(i,j)}$.
By computing $\mL^2$, we find all paths $i - j - k$ such that $i > j > k$ because $\mL$ only contains nonzeros $\mL(i,j)$ where $i> j$ and multiplying it with itself discovers only the path with decreasing vertex id. Computing the Hadamard product of $\mL^2$ with $\mL$ itself ensures the triangle is discovered and only once. Figure~\ref{fig:tc} shows an example.

It might seem obvious that first computing $\mL^2$ in its entirely before doing an element-wise multiplication with $L$ is inefficient for many cases where the number of wedges (i.e. paths of form $i - j - k$) are significantly larger than the number of actual triangles. Ideally, we only want to count wedges around an existing edge. This is why the Masked-SpGEMM has originally been proposed for the triangle counting problem~\cite{trianglegabb15} and has eventually made its way into the GraphBLAS API where the matrix multiplication function {\tt GrB\_mxm}  accepts a mask. The implementation can choose to either drive the computation by going over the existing nonzeros of $\mL$ and computing dot products $\mL(i,:) \cdot \mL(:,j)$ for each $\mL(i,j) \neq 0$, or it can drive the computation by going over the inputs and compute standard SpGEMM except for the outputting part where entries are masked out before writing the result back to memory. 

SpGEMM can also used to perform bulk graph sampling while training graph neural networks using mini-batch sampling. For coherence, this use case is described within Section~\ref{sec:gnns}.

\subsection{SpGEMM for Database Operations}
Different communities have discovered the connection between matrix multiplication and various database joins~\cite{hu2020parallel, kotlyar1997relational, hutchison2017laradb}.
Recent worst-case optimal join algorithms in particular both rely on and improve upon lower bounds on output-sensitive SpGEMM~\cite{amossen2009faster}. 

The connection intuitively follows once we reinterpret matrices $\mA$ and $\mB$ as database tables with two columns each. Then SpGEMM $\mC = \mA \mB$ over the standard real field is 
just a special case of equi-join where (1) the second dimension of $\mA$ and the first dimension of $\mB$ are the columns over which we perform the join, and (2) the duplicates are summed as opposed
to concatenated. Because SpGEMM on semirings can choose to overload the scalar addition operation to perform different operations on duplicates as opposed to merely summing them, 
a generalized SpGEMM over a semiring can easily be used to perform equi-joins. In real life, tables are not restricted to two columns but luckily sparse matrix data structures are ``templated'' in a way
that would allow storing all non-key values on the dimension that is not being used as the contraction dimension during multiplication. Curiously, applications in the opposite direction where the authors advocated using established database technology for matrix computations have been recently proposed~\cite{luo2018scalable}. 

On the practical side, recent work shows benefits of using SpGEMM for context-free path querying on graph databases~\cite{terekhov2020context}. 

\subsection{SpGEMM for Machine Learning}
The Masked-SpGEMM is shown to be a core kernel in the computation of tree-based models for extreme multi-label ranking/classification
(XMR/XMC) problems when sparse features are supported~\cite{etter2022enterprise}. In the linear case of 
XMR, the model is a tree that represents represent the hierarchical clustering of labels. In this particular instance of Masked-SpGEMM $\mY \langle\mM\rangle \gets \mA \mX$, the $n \times L$ output $\mY$ is the 
ranker activations at a particular rank where $L$ is the number of clusters in a particular layer (depth) of the model tree. $\mX$ is the $d \times L$ sparse weight matrix for the tree model at that layer. The $n \times d$ sparse matrix $\mA$ represents the queries. Because the sparse matrix $\mX$ represent the hierarchical tree structure, there is naturally occurring {\emph structured sparsity} that is taken advantage by the authors to develop a masked chunked SpGEMM for this problem. 

\subsection{SpGEMM for Randomized Algorithms}
\label{sec:spgemmrand}
Randomized algorithms have been recently at the forefront of research in numerical linear algebra~\cite{martinsson2020randomized} and have shown growing importance in several areas of scientific computing~\cite{buluc2021randomized}. 
At its core, many randomized algorithms in numerical linear algebra involve the application of a $d$-by-$m$ sketching matrix $\mS$ to an $m$-by-$k$ data matrix $\mA$. 
When both the sketching matrix $\mS$ and the data matrix $\mA$ are sparse, the sketching operation $\mS \mA$ becomes an SpGEMM. However, the data matrix is often tall-and-skinny ($m < k$) with its
rows representing data samples and its columns representing features. The sketching matrix $\mS$ is by construction short with $d \ll m$, making the output $\mS \mA$ a small and dense $d$-by-$k$ matrix.
In fact, this high density of the output is required to sustain approximation bounds of many sketching operations (e.g., the oblivious subspace embedding property~\cite{nelson2013osnap, meng2013low}). Hence, the second case of SpGEMM with the dense output applies, for which there has been recent work~\cite{sobczyk2022pylspack}. 

In the case of Sparse Johnson-Lindenstrauss Transforms (SJLT) of Kane and Nelson~\cite{kane2014sparser}, $\mS$ has exactly $s$ nonzero entries per column. In fact, each nonzero is a Rademacher random variable (i.e. uniform over $\{-1,1\}$). SJLT can be thought as a generalization of 
the CountSketch with $s{=}1$. We have already described the use of Gram matrix computation $\mA\transpose \mA$ in the context of computational biology. The Gram matrix is also needed in randomized algorithms, in particular in the context of least 
squares problems~\cite{sobczyk2022pylspack}. 

\section{The SpMM case}
\label{sec:spmm}

Multiplying a sparse matrix with a tall-skinny dense matrix (SpMM) is often a generalization of sparse-matrix-vector multiplication (SpMV) to the multiple vectors case. Even though the use of this computational kernel go back much further, we were able to trace the acronym SpMM only back to Rich Vuduc's PhD dissertation~\cite{vuducthesis}. 

\subsection{Scientific Computing}

The original motivation to study SpMM was block iterative solvers such as block Lanzcos~\cite{grimes1994shifted} or block Arnoldi~\cite{sadkane1993block}. This type of use case is still in demand for SpMM, especially with the recent popularity of the more robust block eigensolver Locally Optimal Block Preconditioned Conjugate Gradient (LOBPCG)~\cite{knyazev2001toward}. Furthermore, the previously described AMG prolongation and restriction operations can also be described in terms of more regular SpMM operations in lieu of SpGEMM operations at the cost of computing a preprocessing step. This method by McCourt et al.~\cite{mccourt2015sparse} uses graph coloring to find structurally orthogonal columns that can be compressed into one denser column, thus potentially converting one of the matrices into a dense matrix and allowing SpMM computations on them.

Low-rank matrix approximation methods, such as non-negative matrix factorization (NMF) are commonly used when analyzing large and noisy datasets. Such methods have been historically used in recommender systems and have enjoyed applications in various science domains for tasks such as data imputation and dimensionality reduction. The alternating directions approach is a highly popular method for low-rank matrix approximation where we alternate fixing one factor matrix and solving for the other. When the input is sparse, either due to missing entries or actual zeros, the most time consuming part of solution methods based on alternating directions become SpMM. An example of this use case can be found in the MPI-FAUN framework~\cite{kannan2017mpi}.

 \begin{figure}[thb]
  \centering
 \includegraphics[width=0.9\linewidth]{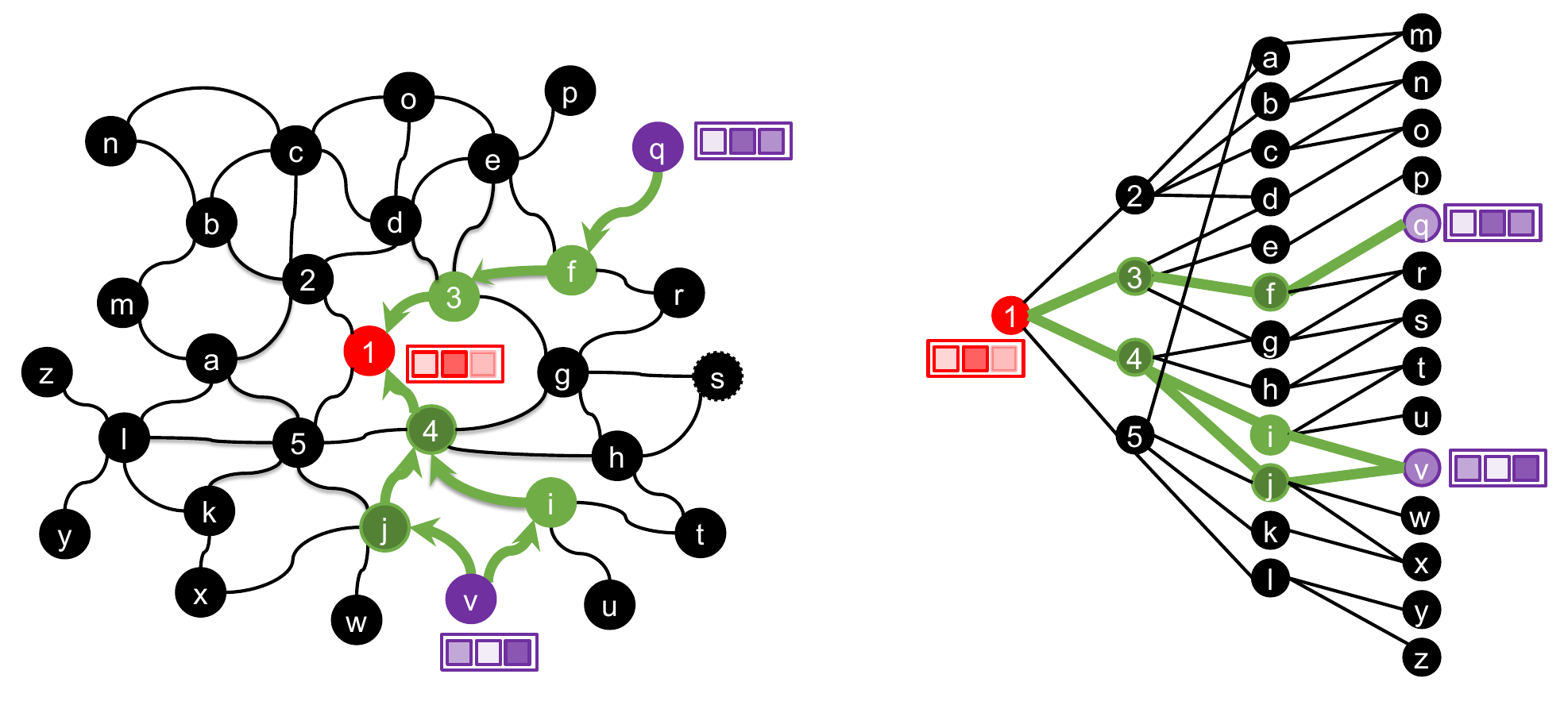}
 \caption{Illustration of the information flow in a Graph Neural Network (GNN). On the left is the graph in its natural form. The features (the shaded boxes) of vertices v and q are aggregated at vertex 1 through intermediate (green) vertices and edges. Features of other nodes are not shown but are also propagated. During training, the error is backpropagated in the opposite direction in the neural network, where each layer of the neural network propagates one hop of information.} 
 \label{fig:gnn}
 \end{figure}

\subsection{Graph Neural Networks}
\label{sec:gnns}
Unlike SpGEMM, an overwhelming majority the known use cases of SpMM is concerned with floating-point arithmetic and does not require an arbitrary semiring. 
Graph Neural Networks (GNNs)~\cite{scarselli2008graph}, highly-successful 
deep learning machines for problems with discrete connected structure, are notable exceptions. 
In addition to the typical trainable matrices that provides an opportunity to mix features across different neural network layers, GNN also allows information to propagate through the network. This way, nodes can learn from their $l$-hop neighbors in an $(l+1)$-layer GNN. This is illustrated in Figure~\ref{fig:gnn}.
SpMM is used in both forward and backward propagation of full-batch GNN training where features are propagated in the graph. Since each object has a k-length feature vector, the propagation is done with an SpMM as opposed to an SpMV. Full expressive power of GNNs can often only be achieved when the functions that (1) aggregate neighbor features and (2) combine that with the existing features of the node itself are flexible~\cite{xu2018how}.
These aggregate and combine operators can be passed as custom add and multiply operations, as allowed by NVIDIA cuSPARSE's \texttt{cusparseSpMMOp()} function or as custom semirings to the SuiteSparse:GraphBLAS.

Training a GNN with full-batch gradient descent has a few downsides: (1) it is known to converge slower to desired accuracy compared to stochastic mini-batch training (SGD), and (2) it uses much more memory than mini-batch training. 

\begin{figure}[thb]
 \centering
 \subcaptionbox{Black nodes are the mini-batch of source vertices chosen for that training step. Vertices with faded shade are not sampled. The rest of the graph is not drawn.}[.52\textwidth]
 {%
 \includegraphics[width=.85\linewidth]{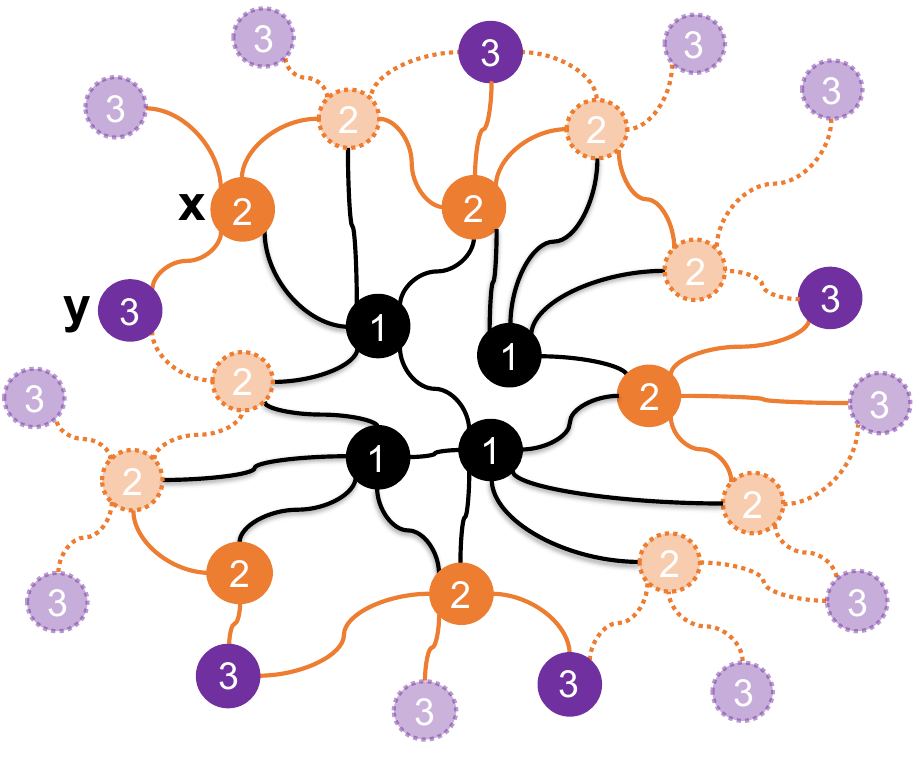}%
 }
 \hfill
 \subcaptionbox{Layer-dependent sampled bipartite graph for the $2$nd layer of GCN. Vertices in $\mathsf{adj}(S_2)$ are highlighted with green halos around them on the right side of the bipartite graph. 
}[.45\textwidth]
 {%
 \includegraphics[width=.9\linewidth]{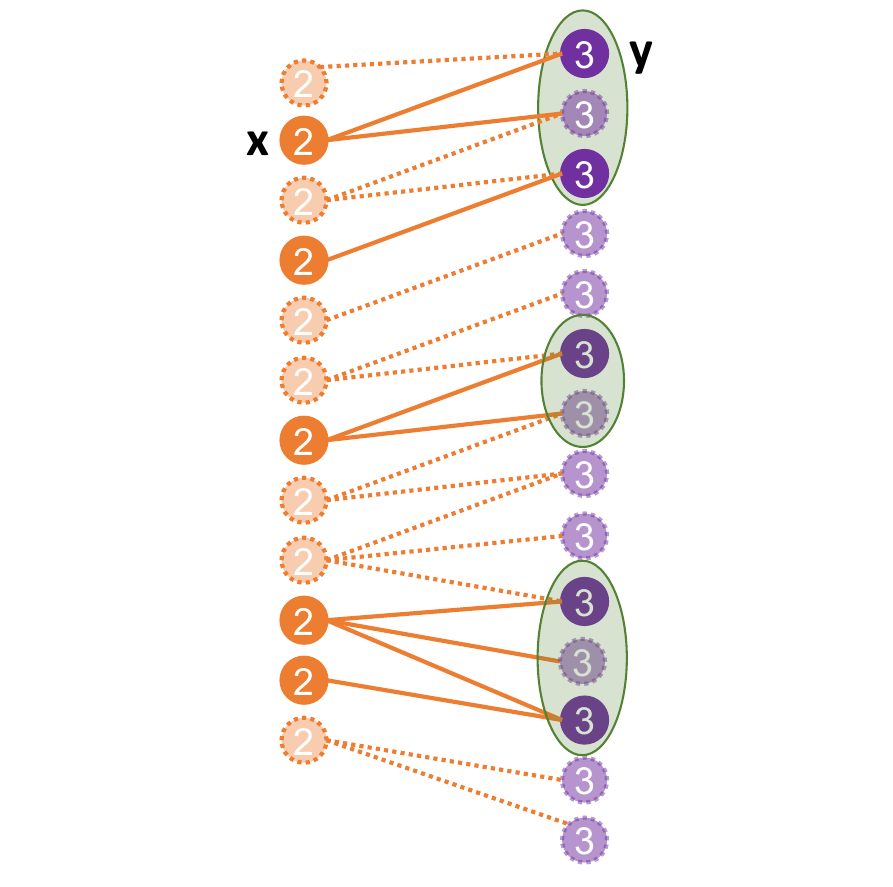}%
 }
 \caption{Layer-dependent sampling illustrated. We sample $5$ vertices at each layer.  Two nodes from layers $2$ and $3$ are marked with $x$ and $y$ as points of reference, respective. The vertices are then ordered within each layer of the bipartite graph clockwise from those reference vertices. 
}
 \label{fig:ladiessampling}
 \end{figure}

The solution is to sample the input graph intelligently during the training process. However, mini-batch training of GCNs is non-trivial because each sample in graph learning is dependent on other samples through edges. Unlike the CNN case where samples are independent, arbitrarily choosing a subset of nodes to mini-batch train a GNN touches almost all the graph, nullifying the benefits of mini-batch training. Various forms of graph sampling methods have been developed to address this issue, but their scalability has been severely limited with inter-processor communication. There are fundamentally three kinds of sampling methods: \emph{node-wise sampling}, \emph{layer-wise sampling}, and \emph{subgraph sampling}. An illustration of layer-dependent neighborhood sampling~\cite{NEURIPS2019LADIES} is demonstrated in Figure~\ref{fig:ladiessampling}. In particular, each layer requires constructing a potentially rectangular submatrix that will serve as the Laplacian matrix for that layer of GCN training. This is computed via $\tilde{\mA} = \mathbf{P}^l \mA \mathbf{D} (\mathbf{P}^{(l+1)})\transpose$. Here, $\tilde{\mA}$ is an $\lvert S_l \rvert \times \lvert S_{l+1} \rvert$ sparse matrix. $\mathbf{P}^l$  and $\mathbf{P}^{l+1}$ are row selection matrices which have one nonzero in each row, and $\mathbf{D}$ is a diagonal matrix. The set of vertices that were sampled in layer $l$ (denoted by $(S_l)$) determine the set of vertices that are eligible for sampling for layer $(l+1)$. In other words, vertices in layer $(l+1)$, denoted by $(S_{l+1})$, are sampled from the aggregated neighborhood of $S_l$, which is denoted by $\mathsf{adj}(S_l)$.

 The operation $\tilde{\mA} = \mathbf{P}^l \mA \mathbf{D} (\mathbf{P}^{(l+1)})\transpose$ can in principle be implemented in parallel using submatrix extraction algorithms based on triple sparse matrix products~\cite{gemmexp}. However, those algorithms' runtimes are proportional to the input graph size. Given that one needs to perform this submatrix extraction repeatedly during each training iteration, \emph{bulk} versions of this operation that extracts multiple sampled submatrices efficiently at once are developed~\cite{tripathy2024distributed}. In these bulk mini-batch sampling algorithms, the calculation of sampling probabilities of vertices is also mapped to an SpGEMM because the sampling probability matrix can be constructed by random walks from multiple groups of vertices, where each group representing a different mini-batch. 

GNNs can also be used to perform deep learning on point cloud data. Point cloud neural networks (PCNNs) are specialized models designed to process and analyze 3D point cloud data -- unordered sets of points representing objects or environments captured by sensors like LiDAR and RGB-D cameras, with applications in robotics, autonomous driving, and AR/VR. Modern PCNNs leverage GNN architectures to model spatial relationships among points. In this pipeline, a k-nearest neighbor graph is built from the input point cloud data, which is then fed to a GNN to make predictions on edges, nodes, or whole graphs. Equivariance, which is the property that the neural network's output transforms (e.g., rotates, translates, reflects) the same way the input transforms, can improve the data efficiency of PCNNs. Equivariant GNNs are particularly effective in quantum chemistry applications, where they can predict molecular properties such as energies and atomic forces for ensembles of molecules~\cite{nequip_2022, Batatia2022mace}.

GNNs have demonstrated significant success in symmetry-preserving equivariant deep learning applied to point clouds -- a finding independently confirmed by several research groups~\cite{cg_nets, thomas_tensor_2018, steerable_cnns}. Specialized algorithms can integrate node and edge features to enforce equivariance with respect to Euclidean transformations of the input space. To achieve this, these methods use a variant of graph convolution whose memory access pattern mirrors that of SpMM. Specifically, a custom operator known as the Clebsch-Gordan (CG) tensor product combines node and edge features before aggregating information across each node's neighborhood. 
Several tools such as OpenEquivariance~\cite{openequivariance}, cuEquivariance~\cite{cuequivariance}, and FlashTP~\cite{lee2025flashtp}, have further optimized rotation-equivariant GNNs by leveraging the SpMM memory access pattern. They fuse the node-edge interaction step directly into the graph convolution, enabling highly efficient computational kernels that avoid the need to explicitly materialize intermediate results, significantly accelerating inference and training.

\subsection{Sparsity in Neural Networks}
The model weights in deep neural networks (DNNs) such as convolutional neural networks or transformers can be sparsified before, during, or after training, using a variety of techniques such as structured pruning. A summary of these sparsification techniques circa 2020 can be found in the excellent survey by Hoefler et al.~\cite{hoefler2021sparsity}. It is challenging to extract performance benefits from this enforced sparsity because pushing the sparsity beyond a certain point is known to degrade accuracy~\cite{gale2019state}. The threshold of sparsity where the accuracy loss becomes problematic is dependent on the deep learning model and is rapidly changing. This \emph{twilight zone of sparsity} in DNNs is often quantified and characterized with percentages of potential entries that are absent. By contrast, the commonly observed sparsity in other applications is characterized by each row or column of a matrix having only a constant or $O(\lg{n})$ nonzeros. This fundamental difference between DNNs and other uses of sparsity happens because the default state is dense in DNNs whereas in almost all other applications described in this paper, the data was already sparsely connected and no explicit sparsification was needed. When the matrix storing the model weights is sparsified, the workhorse of inference  becomes SpMM when mini-batch stochastic gradient descent is employed~\cite{galesc20}. The column dimension of the dense matrix is the mini-batch size, which can vary from 64 to 32K~\cite{you2017scaling}. 
As previously mentioned, this sparse matrix in the twilight zone of sparsity requires non-traditional data structures to extract performance good enough to compete with the heavily optimized dense BLAS operations that take advantage of latest hardware features such as tensor cores. 

Sparsifying weights of a neural network, known as the model sparsity, generally creates fixed sparsity pattern that can be reused in inference repeatedly, justifying expensive transformations to improve inference performance.  Other common forms of sparsity in deep learning are ephemeral, meaning that their form dynamically changes per sample. One example is the sparsification of gradients during optimization using techniques such as keeping the top-k entries~\cite{alistarh2018convergence}. A more commonly observed case of sparsity emerges in activations, especially with certain activation functions such as ReLU and its variants. Recent work found that $90-95\%$ of activations in modern transformers are zero if we look beyond the first layer~\cite{li2022lazy}. 

Finally, SpMM shows up in what is known as a deep learning recommendation model (DLRM)~\cite{naumov2019deep}. The data is categorical and each category is embedded in a $d$-dimensional space, which is concretely represented as a $d$-length vector. The embedding table stores the vector embeddings of all categories, which is a tall-skinny dense matrix $\mathbf{V}$ of dimensions $n \times d$ where $n$ is the number of categories. Given that an item often belongs to multiple categories, a single item lookup requires selecting a (small) subset of rows and computing a weighted average of the embeddings of the categories involved. This can be written as $\mathbf{a_i} \mathbf{V}$ where  $a_{ij} \neq 0$ for if the $i$th item belongs to $j$ category, and $a_{ij} = 0$ otherwise. A recommender system will often batch the lookups for multiple items, say $b$ of them. Hence the resulting operation becomes $ \mathbf{W} = \mA \mathbf{V}$ where the $i$th row of $\mA$ is $\mathbf{a_i}$~\cite{naumov2019dimensionality}. In contrast to previous use cases, the sparse matrix $\mA \in \mathbb{R}^{b \times n}$ can be short and fat (i.e., $b \ll n$), especially if the batch size is small. 

\subsection{Randomized Algorithms}

Earlier in Section~\ref{sec:spgemmrand}, we covered the case of sketching a sparse data matrix with a sparse sketching matrix using SpGEMM. 
When the data matrix is dense, sparse sketches can be applied 
to it using SpMM with the same caveats. 
In theory, another potential use of SpMM for sketching would be to apply a dense sketch to a sparse data matrix. The trouble with dense sketch matrices is the enormous cost of storing them. Instead,
a better approach would be to just regenerate the elements of the dense sketching matrix on the fly as needed~\cite{liang2024fast}. Murray et al.~\cite{murray2023randomized} provides a comprehensive treatment of computational aspects of randomized linear algebra 

Having briefly mentioned the use of SpMM in sketching, we now turn our attention to randomized algorithms outside numerical linear algebra.
Consider the following problem that is fundamental for genome assembly using the latest sequencing technology (i.e., single-molecule sequencing (SMS)). The output of sequencing instruments is a set of \emph{reads}, which themselves are generally orders of magnitude shorter than the genome that needs to be assembled. The first step of the assembly process using the SMS technology is to find the overlaps between noisy reads. The brute force all-to-all comparison of read pairs is infeasible partly because comparing a pair of reads is expensive: quadratic in the length of those reads. Consequently, we would like to avoid as many unnecessary read-to-read comparisons as possible. An effective method is to use $k$-mers, very short nucleotide sequences of length $k$. It can be shown that for an appropriately chosen length $k$, which depends on the error rate and the average read length, two reads are very unlikely to overlap if they do not share any pair of common $k$-mers.

One can count $k$-mers in the whole dataset and identify unique $k$-mers that are present in each read. The absence and presence of $k$-mers in a read serve as the feature vector of that read. We can then compute the read-to-read similarity using the number of shared $k$-mers. The problem is that there are $4^k$ possible distinct $k$-mers using the DNA alphabet and an overwhelming majority will be encountered after sufficient sequencing depth, given the high error rates of SMS technology. 

\begin{figure}
\begin{center}
\includegraphics[scale=0.45]{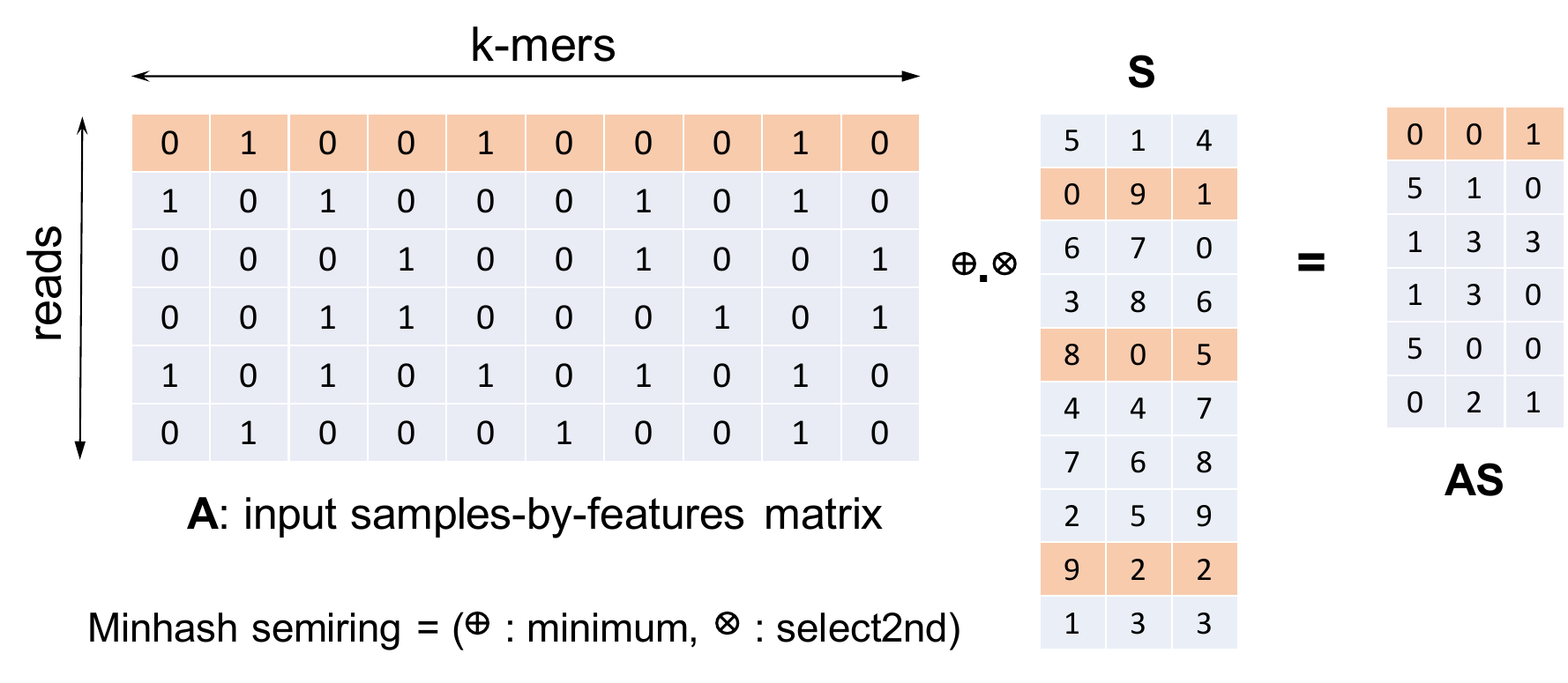}
\caption{Application of minhashing on a dataset with each row (sample) being a DNA read and each column (features) being its $k$-mers. We perform three min-wise independent permutations on the columns of this data to produce a sketch that preserves the distances between samples. Minhashing here is implemented exactly without explicit column permutations. Instead, a sparse-matrix times dense-matrix over an algebraic semiring is used. The multiplication operation of the semiring returns the second operand and the addition returns the minimum. The rows accessed during the formation of the first row of the output are highlighted in orange.
}
\label{fig:minhash}
\end{center}
\end{figure}

A randomized approach called minhashing~\cite{broder2000min} is commonly used to reduce the large feature space drastically. This approach is illustrated in Figure~\ref{fig:minhash} where the input data is represented as a $\lvert \textit{reads} \rvert \times \lvert \textit{k-mers} \rvert$ sparse matrix $\matrix{A}$. The minhashing process is commonly described as applying $r$ independent random permutations to columns of $\mA$. These random permutations are stacked as columns of the $\matrix{S}$ matrix in our example. After the $i$th permutation is applied, the column index of the first nonzero in each row (sample) serves as the $i$th component of that sample's signature. Existing work shy away from this exact application of minhashing due to high costs of explicit permutation of the data matrix and instead resort to approximations. However, a closer look reveals that an explicit permutation is not computationally necessary for exact minhashing. In our example, the first row has only three nonzeros with column indices $\{1,4,8\}$ and the first permutation sends those nonzeros to $\{0,8,9\}$th locations. Consequently, the first component of the first row's signature is $0$ because that is the minimum column index of the nonzeros in that row post permutation.

As seen in Figure~\ref{fig:minhash}, minhashing can be implemented using SpMM (or SpDM$^3$ depending on the aspect ratio of the dense matrix) over an algebraic semiring. Notice that each column of the matrix $\matrix{S}$ is a permutation of $\{0,1,\ldots, 9\}$. The similarities between rows of the output matrix $\mA \matrix{S}$ approximately preserve the similarities between rows of the input matrix $\matrix{A}$ but the similarity computation can now be carried out at the fraction of the cost. GraphBLAS provides flexible support for semirings that will be used as a tool to implement minhashing directly using SpMM. We note that minhashing is designed and widely used within the information retrieval community; for example to find near-duplicate documents in the web. In that context, the boolean sparse matrix is the documents-by-shingles matrix where a shingles is the set of strings of length $k$ that appear in the document, similar to a $k$-mer.

\begin{figure}
\begin{center}
\includegraphics[scale=0.5]{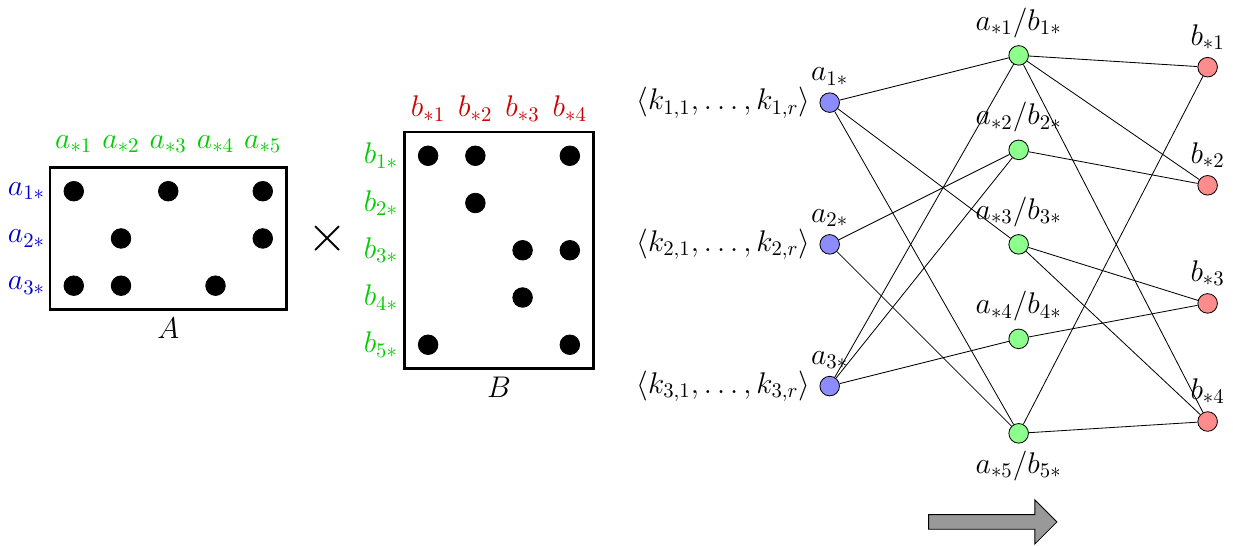}
\caption{Illustration of the layered graph concept used in Cohen's algorithm for nonzero estimation in SpGEMM. Each propagation of the random values (the vector of keys $k_{i,1}
\ldots, k_{i,r}$) can be implemented via a SpMM on a (select2nd,min)-semiring. The figure is replicated from earlier work~\cite{selvitopi2020optimizing}.
}
\label{fig:nonzeroestimation}
\end{center}
\end{figure}

Another application of SpMM in the land of combinatorial randomized algorithms is going to bring us full circle because the algorithm predicts the nonzero structure of SpGEMM. In particular, it predicts row (or column)
nonzero counts of the output matrix. It is also trivially applicable to matrix chain products, hence making it especially useful for graph contraction with triple products we described in earlier sections.  
Cohen's algorithm~\cite{cohen1998structure} models the sparse matrices and their product using a layered graph. Each level of the graph corresponds to the rows or columns involved in the matrix multiplication. Nonzero entries are represented as edges between nodes in adjacent levels, capturing the potential pathways that lead to nonzero entries in the product matrix. This is illustrated in Figure~\ref{fig:nonzeroestimation}. 

Let us assume for the sake of illustration, we are estimating column nonzero counts for $\mC = \mA \mB$. The algorithm randomly generates probability values for all source vertices (which represent rows of the first matrix $\mA$). These values are propagated through this graph (analogous to breadth-first search from all source vertices), where each internal or destination vertex aggregates the minimum value it received. The number of distinct source vertices that can reach each destination vertex (representing columns of the second matrix $\mB$) can be used as a statistical estimator for that column's nonzero counts, after necessary transformations. To understand why this approach successfully estimates the number of nonzeros in each column instead of flops, consider two scenarios. 

\begin{enumerate}[leftmargin=10pt]
\item If there are two intermediate products $a_{ik} b_{kj}$ and $a_{il} b_{lj}$ contributing to the same nonzero $c_{ij}$, then there will be two distinct paths between the source vertex $a_{i*}$ that represents the $i$th row of $\mA$ and the destination vertex $b_{*j}$ that represents the $i$th column of $\mB$. Since only one unique random number (that is, $a_{i*}$'s) reaches $b_{*j}$ (albeit from distinct paths), this does not effect the statistical estimate compared to having only one intermediate product leading to $c_{ij}$. This is exactly what we want. 
\item When two distinct source vertices $a_{k*}$ and $a_{l*}$ can reach the same column $b_{*j}$, this means that there are two distinct nonzeros $c_{kj}$ and $c_{lj}$ being generated in the output, regardless of the intermediate vertex their paths pass through. Because we keep the minimum when these paths collide at $b_{*j}$, it does effect the statistical estimate (in favor of increased nonzero count in the $j$th column of $\mC$). This is also the desired outcome. 
\end{enumerate}

To implement this algorithm efficiently without explicitly constructing an auxiliary graph, one can generate a random vector of keys $\mathbf{k}$, drawn from an exponential distribution. This vector would have length of length $M$, which is the number of rows of $\mA$. We then perform two sparse matrix-dense vector multiplications $\mathbf{z}= \mB^T (\mA^T \mathbf{k})$ on the (select2nd,min)-semiring. 

Cohen's estimator (and similar randomized sampling approaches) needs many independent random test vectors to decrease the variance. We could do multiple SpMV runs but that is inefficient. 
Packing the vectors $\mathbf{k}$ as columns of a dense matrix, running a single SpMM, and post-processing columns achieves identical mathematical results to multiple SpMV runs, but is dramatically faster in practice.

\section{The SDDMM case}
\label{sec:spdmmm}

 \begin{figure}[thb]
 \centering
 \subcaptionbox{Self-attention mechanism for Vertex 1 for a GAT with feature embedding dimension of 3. Different shades signify different values or magnitudes of the features.}[.55\textwidth]{%
 \includegraphics[width=.75\linewidth]{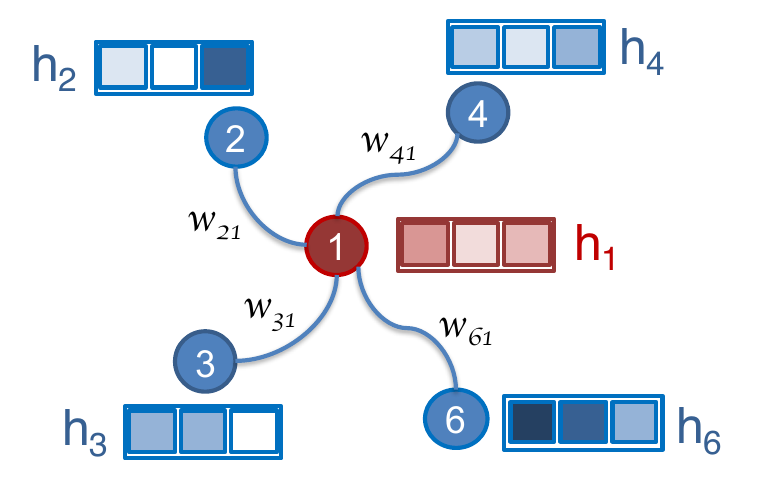}%
 }\hfill
 \subcaptionbox{A few examples for computing the attention weights using the dot product.}[.40\textwidth]{%
 \includegraphics[width=.75\linewidth]{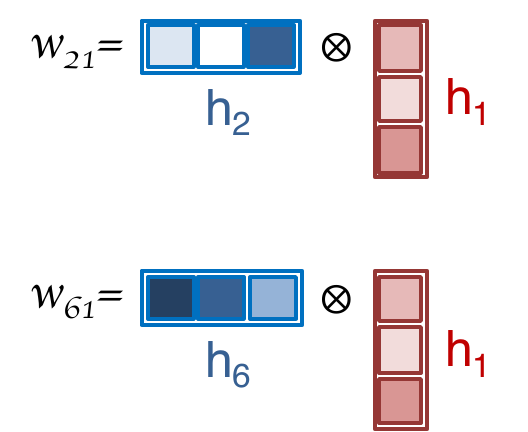}%
 }
  \subcaptionbox{Computing full-batch attention using Masked GEMM, aka sampled dense-dense matrix multiplication (SDDMM). $\otimes$ is multiplication on a semiring, $\circ$ is the Hadamard product.}[.85\textwidth]{%
 \includegraphics[width=.85\linewidth]{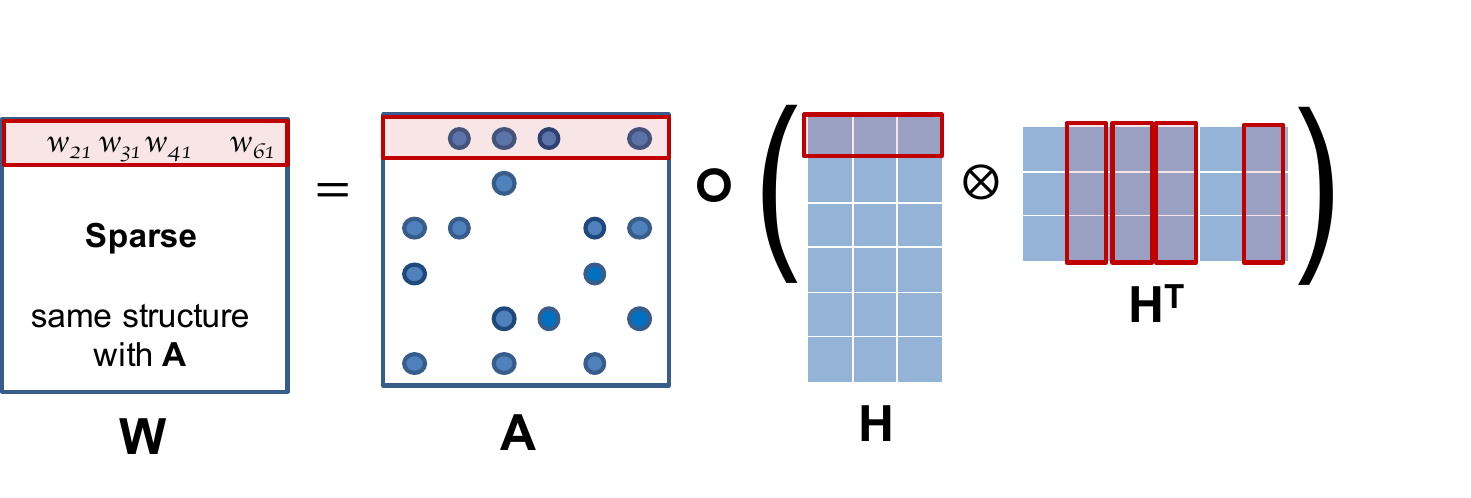}%
 }
 \caption{Graph Attention Mechanism and its efficient batched implementation using Masked GEMM (also known as SDDMM). $\mathbf{W}$ can then be uses as input to a (weighted) graph convolution.
}
 \label{fig:gat}
 \end{figure}

GCNs require the input to explicitly state the strengths of the connections between graph vertices (i.e., edges in the graph) and these connections typically do not change during training. Graph Attention Networks (GATs)~\cite{VelickovicCCRLB18} are more flexible as they also learn the strength of connections between vertices during training. This is achieved by the self-attention mechanism, which learns the weight of each edge using the embedding features of its vertex endpoints. The process is illustrated in Figure~\ref{fig:gat} for the case where attention is a dot product. Multi-head attention further enhances the power of GATs by learning multiple weights simultaneously. 

In transformer models, the attention score is computed via the formula

$$ \mathit{Attention}(Q,K,V)=\mathit{softmax} \left( \frac{Q K^T}{\sqrt{d_k}} \right) V, $$

where $Q$, $K$, $V$ are matrices representing query, key, value embeddings and $d_k$ is the dimensionality of the keys.  It has been observed that one can 
compute a subset of the $O(n^2)$ attention scores for a sequence of length $n$, and still achieve good accuracy. The sparsity pattern of the attention matrix
can be statically or dynamically determined. Static sparse attention chooses which tokens to attend using their relative positioning (e.g., sliding window), whereas
dynamic sparse attention chooses tokens to attend by computing an importance score. Several influential sparse attention methods have been proposed for long context attention, 
including Longformer~\cite{beltagy2020longformer}, Big Bird~\cite{zaheer2020big}, 
and sparse transformers~\cite{child2019generating}. These methods can reduce the quadratic cost to compute attention scores to linear in the sequence length. 

Algorithmically, if we consider the sparsification method as a sparse mask $M$, the attention computation becomes 
$$ \mathit{Attention}(Q,K,V)=\mathit{sparse\_softmax} \left( \frac{M \odot Q K^T}{\sqrt{d_k}} \right) V.$$

Here, $A = M \odot Q K^T$ (or $A \langle M \rangle \gets Q K^T$) is an SDDMM operation.
Since the temporary output matrix $A$ is sparse, $B = \mathit{sparse\_softmax} (A / \sqrt{d_k})$ is also sparse. Hence the subsequent multiplication of $B$ with $V$ is now an SpMM operation.
This is recognized by Liu et al, who went on to show that using cuSPARSE for sparse attention can beat dense attention for $>80\%$ sparsity in full precision~\cite{liu2022dynamic}. However, most deep learning pipelines
are in low precision and it requires significant engineering effort to beat low-precision dense kernels. DeepSeek recently demonstrated~\cite{yuan2025native} that their natively implemented sparse attention, which uses
a combination of dynamic and static attention, can significantly outperform state-of-the-art attention implementations, such as FlashAttention-2. In particular, native sparse
attention was $9\times$ faster in forward propagation and $6\times$ faster in backpropagation for 64K context length in training. The improvements were more pronounced as the context length gets larger.

\section{The SpDM$^3$ case}
We have already covered the multiplication of sparse matrix with a tall-skinny dense matrix (SpMM). However, there are cases where the dense matrix is also square and hence is drastically larger than the sparse matrix (in terms of its nonzero count). This use case, which we call SpDM$^3$ following earlier work~\cite{koanantakool2016communication}, presents sufficiently distinct challenges for us to cover it separately from SpMM.
Admittedly, the applications of SpDM$^3$ are less common than all the other primitives we cover in this paper.

Sparse inverse covariance matrix estimation is a statistical method used to estimate the inverse of a covariance matrix (also called the precision matrix) where many entries are constrained to be zero. This technique is widely used in high-dimensional data analysis, such as graphical models, because a sparse precision matrix corresponds to conditional independence structure among variables. There are many popular methods for estimating the precision matrix, including those based on maximum likelihood e.g. (BigQUIC~\cite{hsieh2013big}, SQUIC~\cite{bollhofer2019large}), as well as those based on pseudolikelihood (e.g., CONCORD~\cite{khare2015convex}). All methods are $l1$-regularized to enforce sparsity in the solution. CONCORD minimizes the following objective function:

$$
\min_{S} \; \mathrm{tr}(Q S^2) - 2 \log \det(\mathrm{diag}(S)) + \lambda \sum_{i \neq j} |s_{ij}|,
$$

\noindent
where $\mathbf{Q}$ is the sample covariance matrix, which is often dense. It can be computed from the data matrix $\mathbf{X}$ via $\mathbf{Q}= \frac{1}{n} \mathbf{X} \mathbf{X}^T$.
Sparse inverse covariance matrix $\mathbf{S}$ is the object we want to estimate.
The CONCORD method and its variants~\cite{oh2014optimization}, including a distributed-memory parallelization HP-CONCORD~\cite{koanantakool2018communication}, needs to 
compute $\mathbf{S}^k \mathbf{Q}$ on each proximal gradient and line search iteration where $\mathbf{S}^k$ is the $k$th estimate for the precision matrix. Each iterate $\mathbf{S}^k$ is sparse, resulting in the
SpDM$^3$ function calls. 

All pairs shortest paths (APSP) is a classical graph problem where the goal is to compute the shortest paths amongst all pairs of vertices. The famous Floyd-Warshall algorithm solves this problem in optimal $O(n^3)$ time for a dense graph with $n$ vertices. Using the $(min,+)$ semiring, the distance only version of APSP can be calculated via matrix multiplications, either via the path-doubling method~\cite{seidel1995all} or by the recursive Kleene's formulation~\cite{d2007r}. The shortest paths themselves can also be computed on alternate semirings using path algebra at the same cost~\cite{tarjan1981unified, hofner2012dijkstra}. 
The recursive Kleene's algorithm is also optimal for dense graphs, and is especially well suited for high-performance computers such as GPUs~\cite{bulucc2010solving} and distributed-memory clusters~\cite{solomonik2013minimizing}. 
The classical path doubling method, while offering tremendous parallelism, has an additional $log(n)$ multiplicative cost that makes is less competitive in practice. 

While all the aforementioned algorithms for APSP are presented for dense graphs, they actually are fairly competitive for sparse graphs as well due to the fact that the output is often dense. 
However, the intermediate calculations can in theory take advantage of sparsity. Tiskin was the first to recognize and exploit this observation~\cite{tiskin2001all}. He realized that any shortest path with up to 
$2k$ edges can be represented as either a shortest path with up to $k$ edges (already computed), or a path that consists of a shortest path of exactly $k$ edges (from some vertex $u$  to some $z$), followed by another shortest path of up to $k$ edges (from $z$ to $v$). Hence, the main distance relaxation formula becomes:

$$
D^{(2k)} = D^{(k)} \oplus \left( D^{[k]} \otimes D^{(k)} \right)
$$

\noindent
where $\oplus$ denotes element-wise minimum, and $\otimes$ denotes min-plus matrix multiplication. $D^{(k)}$
  is the matrix of shortest path distances for paths of up to 
$k$ edges. $D^{[k]}$
  is the matrix of shortest path distances for paths of exactly 
$k$ edges. $D^{[k]}$ matrix is typically sparse, allowing us to utilize SpDM$^3$ kernels during each iteration, as described by Solomonik and Hoefler~\cite{solomonik2015sparse}.
Sparsity in $D^{[k]}$ emerges because as 
$k$ increases, the number of shortest paths of exactly 
$k$ edges becomes much less than the total number of all (up to 
$k$-edge) paths. 

\section{Conclusions and Future Outlook}
\label{sec:conclusions}

A growing number of data analytics and scientific computations can be carried out efficiency with the sparse matrix abstraction.
In this paper, we reviewed the sparse matmul operation and its applications.
We identified the following features of sparse matmul as most important for new algorithm and architectural studies. 

For both SpGEMM and SpMM, the implementations should be capable running on a larger set of semirings than just floating-point arithmetic. This need has been well documented for graph algorithms but the recent success of graph learning techniques such as GNNs and graph kernels make it a more pressing need. For example, FeatGraph~\cite{featgraph} relied on the TVM compiler~\cite{chen2018tvm} to compile user-defined functions for GNNs. We expect more refinements and generalization of this approach thanks to the increasing sophistication of sparse tensor algebra compilers~\cite{senanayake2020sparse}.   

Sparse deep learning operates on a delicate zone where the sparsity cannot be pushed too much without losing accuracy but the potential memory and performance benefits are too important for the developers to ignore. Consequently, new data structures and algorithms should be developed for targeting this \emph{twilight zone of sparsity} that is arising in traditional deep learning.  In particular, the index storage overhead is worse in deep learning due to widespread use of low-precision formats~\cite{micikevicius2022fp8}. This necessitates increased use of index compression~\cite{willcock2006accelerating} and hierarchical blocking schemes that limit the addressable index space~\cite{bulucc2009parallel}. We are encouraged by the recently proposed sparse matrix formats that address sparse deep learning requirements~\cite{lin2025toward}. One clear advantage of sparsity research in deep learning is the high computational burden of deep learning training and inference. Many sparse matrix reordering and blocking methods designed for other use cases failed to gain traction as the cost for reordering and/or blocking was not amortized by the gains from the application use case. The cost of deep learning can tilt the scales towards more expensive reordering and/or blocking algorithms. On the negative side, one can argue that sparse computations so far have been the losers of the hardware lottery~\cite{hooker2021hardware}, except for a handful of attempts such as the introduction of sparse tensor cores by NVIDIA~\cite{mishra2021accelerating}. The hardware landscape can change, however, with more promising results such as DeepSeek's native sparse attention. 

Single-node CPU and GPU developers should think about the use case of their kernels within a larger-scale distributed code and henceforth worry about hypersparsity as one of their benchmarks. Implementations that are overly optimized as single-node CPU or single GPU kernels tend to perform poorly when used as a subroutine in a large distributed code. 

\section*{Acknowledgments}
This work is supported by the Advanced Scientific Computing Research (ASCR) Program of the
Department of Energy Office of Science under contract No. DE-AC02- 05CH11231.

We thank John D. Owens for kick-starting the writing of this article by asking the question of real-life applications of SpGEMM in an email. 
We thank Vivek Bharadwaj for his help with the SpMM's use within equivariant neural networks. 

\bibliographystyle{siamplain}
\bibliography{references}

\end{document}